\documentclass{amsart}

\usepackage[latin1]{inputenc}
\usepackage[english]{babel}
\usepackage{indentfirst}
\usepackage{amssymb}
\usepackage{amsthm}
\usepackage{xcolor}
\usepackage[all]{xy}

\newcommand{\N}{\mathbb{N}}

\newcommand{\sub}{\subseteq}

\newcommand{\ca}{\operatorname{ca}}

\def\epsilon{\varepsilon}

\newtheorem{theo}{Theorem}[section]
\newtheorem{lem}[theo]{Lemma}
\newtheorem{pro}[theo]{Proposition}
\newtheorem{cor}[theo]{Corollary}

\numberwithin{equation}{section}

\title{Ranges of vector measures and copies of $c_0$}

\author{Jos\'{e} Rodr\'{i}guez}

\address{Dpto. de Matem\'{a}ticas\\E.T.S. de Ingenier\'{i}a Agron\'{o}mica y de Montes y Biotecnolog\'{i}a\\
Universidad de Castilla-La Mancha\\ 02071 Albacete\\ Spain} 
\email{jose.rodriguezruiz@uclm.es}



\thanks{The research was supported by grants PID2021-122126NB-C32 
(funded by MCIN/AEI/10.13039/501100011033 and ``ERDF A way of making Europe'', EU) and 
21955/PI/22 (funded by {\em Fundaci\'on S\'eneca - ACyT Regi\'{o}n de Murcia}).}

\begin{document}

\begin{abstract}
Let $X$ be a Banach space such that every weakly compact subset of~$X$ is contained in the range of an $X$-valued countably additive
vector measure defined on a $\sigma$-algebra. We prove that $X$ contains no isomorphic copy of~$c_0$, thereby 
answering a question posed by R.~Anantharaman and J.~Diestel [Comment. Math. Prace Mat. 30 (1991), no. 2, 221--235].
In particular, every Banach lattice with this property is superreflexive. 
\end{abstract}

\maketitle

\section{Introduction}

By a {\em vector measure} we mean a countably additive map defined on a $\sigma$-algebra and taking values in a Banach space.
Let $X$ be a Banach space, let $\Sigma$ be a $\sigma$-algebra and let $\nu:\Sigma\to X$ be a vector measure. A classical result
of Bartle, Dunford and Schwartz (see, e.g., \cite[p.~14, Corollary~7]{die-uhl-J}) states 
that the range of~$\nu$, that is, the set 
$$
	\mathcal{R}(\nu)=\{\nu(A): \, A\in \Sigma\}
$$
is relatively weakly compact in~$X$. In fact, Diestel and Seifert~\cite{die-sei} proved that this set has a stronger property, namely, the Banach--Saks property:
every sequence in~$\mathcal{R}(\nu)$ admits a subsequence whose arithmetic means are norm convergent.
Considerable effort has been devoted to studying ranges of vector measures, in particular, the 
question of which sequences may occur in ranges of 
vector measures. For detailed information on this topic, see~\cite{sof} and the references therein.

Anantharaman and Diestel~\cite[p.~230]{ana-die}
asked which Banach spaces $X$ have the property, hereafter called {\em property~AD}, that every weakly compact subset of~$X$ 
is contained in the range of an $X$-valued vector measure; more specifically, they asked whether $c_0$ enjoys property~AD. We next comment on the known facts about 
this property, which explain the relevance of this question. 

Clearly, property AD is preserved by isomorphisms and is inherited by complemented subspaces.
By the result of Diestel and Seifert, a Banach space $X$ with property~AD also
enjoys the \emph{weak Banach--Saks property} (that is, every weakly compact subset of~$X$
has the Banach--Saks property or, equivalently, every weakly null sequence in~$X$ admits a subsequence whose arithmetic means are norm null).
For instance, the spaces $c_0$ and $L_1[0,1]$ have the weak Banach--Saks property (see, e.g., \cite[Exercise~3.57]{fab-ultimo} and \cite[p.~112]{die-J}).
 
All finite-dimensional Banach spaces have property~AD. More interestingly, the spaces
$\ell_p$ and $L_p[0,1]$ have property~AD whenever $2\leq p <\infty$. Indeed, on the one hand,
it is known that the closed unit ball $B_X$ of a Banach space~$X$ {\em is} the range of an $X$-valued vector measure if and only if 
its dual $X^*$ (with the dual norm) is isometric to a reflexive subspace of~$L_1(\mu)$ for some probability measure~$\mu$ 
(see \cite[Proposition~A]{ana-die}). On the other hand, this holds true if $X$ is any of the spaces $\ell_p$ and $L_p[0,1]$ for $2\leq p<\infty$,
since $\ell_q$ and $L_q[0,1]$ embed isometrically into~$L_1[0,1]$ for $1<q\leq 2$ 
(see, e.g., \cite[Proposition~11.1.9]{alb-kal}). 

However, for $1\leq p<2$ the space~$\ell_p$ fails property~AD (see \cite[Examples~6 and~8]{ana-die}) and so does
$L_p[0,1]$ (which contains $\ell_p$ as a complemented subspace).

A reflexive Banach space~$X$ has property~AD if and only if 
$B_X$ {\em is contained in} the range of an $X$-valued vector measure. 
It is known that such a space is superreflexive. Indeed, the proof of
\cite[Proposition~A]{ana-die} shows that if there is an
 $X$-valued vector measure whose range contains~$B_X$, then $X^*$ is isomorphic to a reflexive subspace of $L_1(\mu)$ for some 
probability measure~$\mu$. A result of Rosenthal~\cite{ros-J-5} (cf. \cite[Chapter~IV, \S4]{li-que}) then applies to conclude that $X^*$ is isomorphic to 
a subspace of $L_q(\mu')$ for some $1<q\leq 2$ and for some probability measure~$\mu'$. So, $X^*$ is superreflexive and 
therefore the same holds for~$X$ (see, e.g., \cite[Chapter~9]{fab-ultimo}). For more examples of superreflexive Banach spaces whose
closed unit ball is not contained in the range of a vector measure, see~\cite{cas-san}.

Beyond the reflexive case, it was shown in \cite[pp.~229--230]{ana-die} that the separable infinite-dimensional $\mathcal{L}_\infty$ space 
with the Schur property constructed by Bourgain and Delbaen~\cite{bou-del} has property~AD. The argument uses the so-called weakly compact extension property
and can be extended to prove that every $\mathcal{L}_\infty$ space with the Schur property has property~AD. Since $\ell_1$ fails property~AD,
this shows that, in general, property AD is not inherited by subspaces.

In this paper we address the aforementioned question of Anantharaman and Diestel and we show that
{\em $c_0$ fails property~AD} (Theorem~\ref{theo:c0}). Moreover, by refining some arguments of De Kock and~Puglisi~\cite{koc-pug}, 
the example of $c_0$ can be pushed further to prove that {\em every Banach space having property~AD contains no isomorphic copy of~$c_0$} 
(Theorem~\ref{theo:subspace}). As an application, it follows that {\em all Banach lattices with property AD are superreflexive} (Corollary~\ref{cor:lattice}).

\subsection*{Terminology}
We write $\mathbb N=\{1,2,\dots\}$. The cardinality of a set $F$ is denoted by~$|F|$. 
An {\em operator} is a linear continuous map between Banach spaces. The norm of a Banach space~$X$ is denoted by $\|\cdot\|_X$ or simply $\|\cdot\|$.
The linear subspace generated by a set $C \sub X$ is denoted by ${\rm span}(C)$ and its norm closure is denoted by $\overline{{\rm span}}(C)$.
A norm closed linear subspace of~$X$ is called a {\em subspace}.  
A {\em projection} on~$X$ is an operator $P:X \to X$ such that $P\circ P=P$.

\section{The results}\label{section:results}

We identify $c_0^*=\ell_1$ and, for each $n\in \N$, we denote by $e_n^*\in B_{c_0^*}=B_{\ell_1}$ the $n$-th coordinate functional.

The weakly compact subset of~$c_0$ witnessing the failure of property~AD can be described explicitly. 
Let 
\[
 \mathcal S:=\{\emptyset\}\cup
 \{F\subseteq\N: \, F\text{ is finite and }|F|\leq \min F\}
\]
be the {\em Schreier family} and let 
\[
 K_{\mathcal S}:=\{\chi_F: \, F\in\mathcal S\}\subseteq B_{c_0},
\]
where $\chi_F$ denotes the characteristic function of~$F$. 
\begin{lem}\label{lem:KS}
$K_{\mathcal{S}}$ is weakly compact in~$c_0$.
\end{lem}

In the proof of Lemma~\ref{lem:KS} we will use the following folk result:

\begin{lem}\label{lem:weakconvergence}
Let $X$ be a Banach space, let $\Gamma \sub X^*$ be a set such that 
$X^*=\overline{{\rm span}}(\Gamma)$ and let $\sigma(X,\Gamma)$ be the topology on~$X$ of pointwise convergence on~$\Gamma$.
Then the restrictions of the weak topology and $\sigma(X,\Gamma)$ to every 
bounded subset of~$X$ coincide.
\end{lem}

\begin{proof}[Proof of Lemma~\ref{lem:KS}]
Observe first that $K_{\mathcal{S}}$ is a closed (hence, compact) subset of~$\{0,1\}^\mathbb{N}$ with the product topology. Indeed,
if $F \sub \N$ and $F\not\in \mathcal S$, then $F$ contains at least $\min F+1$ elements.
We can take a finite set $F_0 \sub F$ with $|F_0|=\min F+1$ and $\min F_0=\min F$,
and then $\{\chi_G:  \, G \sub \N \text{ and }F_0 \sub G\}$ is an open neighborhood of~$\chi_F$ in $\{0,1\}^\mathbb{N}$ disjoint from~$K_{\mathcal S}$.

Write $\Gamma:=\{e_n^*:n\in \N\} \sub B_{\ell_1}$. 
Clearly,  
the relative topology of~$K_{\mathcal S}$ as a subset of~$\{0,1\}^\N$ coincides with the restriction of $\sigma(c_0,\Gamma)$ to~$K_{\mathcal S}$.
Since $\ell_1=\overline{{\rm span}}(\Gamma)$, the restrictions of the weak topology and $\sigma(c_0,\Gamma)$ to~$K_{\mathcal S}$ coincide
(Lemma~\ref{lem:weakconvergence}). Hence, $K_{\mathcal{S}}$ is weakly compact in~$c_0$.
\end{proof}

\begin{theo}\label{theo:c0}
$K_{\mathcal{S}}$ is not contained in the range of any $c_0$-valued vector measure. Consequently, $c_0$ fails property~AD.
\end{theo}

Our proof of Theorem~\ref{theo:c0} uses Szlenk's theorem~\cite{szl} that $L_1(\mu)$ has the weak Banach--Saks property
for every positive measure~$\mu$ (non necessarily finite). While the result is often stated just for the Lebesgue measure on the unit interval
(see, e.g., \cite[p.~112]{die-J}), the proof works similarly for an arbitrary measure space.
Actually, we will need a version of Szlenk's theorem in spaces of scalar-valued measures, see
Theorem~\ref{theo:Szlenk} below (the proof is included for the sake of completeness).

Given a measurable space $(\Omega,\Sigma)$, we denote by ${\rm ca}(\Sigma)$ the Banach space of all scalar-valued countably additive maps~$\xi$ defined on~$\Sigma$, 
equipped with the total variation norm $\|\xi\|_{{\rm ca}(\Sigma)}:=|\xi|(\Omega)$. The following characterization of its weakly convergent sequences
(see, e.g., \cite[p.~90, Theorem~11]{die-J}) will be useful for us:

\begin{theo}\label{theo:ca}
Let $(\Omega,\Sigma)$ be a measurable space. A sequence $(\xi_n)_{n\in \N}$ in~${\rm ca}(\Sigma)$ converges weakly to $\xi\in {\rm ca}(\Sigma)$
if and only if $\lim_{n\to \infty}\xi_n(A)=\xi(A)$ for every $A\in \Sigma$.
\end{theo}

\begin{theo}\label{theo:Szlenk}
Let $(\Omega,\Sigma)$ be a measurable space. Then
${\rm ca}(\Sigma)$ has the weak Banach--Saks property.
\end{theo}
\begin{proof}
Let $(\xi_{n})_{n\in \N}$ be a weakly null sequence in~${\rm ca}(\Sigma)$. 
Let $\mu$ be a finite positive measure on~$\Sigma$
such that each $\xi_n$ is $\mu$-continuous (that is, $\xi_n(A)=0$ for every $A\in \Sigma$ with $\mu(A)=0$); for instance, one can take
$$
	\mu:=\sum_{\xi_n\neq 0}\frac{1}{2^n |\xi_n|(\Omega)}|\xi_n|
$$
(the series converges in total variation).

For each $n\in \N$, let $f_n \in L_1(\mu)$ be the Radon--Nikod\'{y}m derivative of~$\xi_n$ with respect to~$\mu$. Since
$(\xi_{n})_{n\in \N}$ is weakly null in~${\rm ca}(\Sigma)$, we have
\begin{itemize}
\item[(i)] $\sup_{n\in \N}\|f_n\|_{L_1(\mu)}=\sup_{n\in \N}\|\xi_n\|_{{\rm ca}(\Sigma)}<\infty$;
\item[(ii)] $\lim_{n\to \infty} \int_A f_n\, d\mu=\lim_{n\to \infty} \xi_n(A)=0$ for each $A\in \Sigma$ (Theorem~\ref{theo:ca}).
\end{itemize} 
Write $\Gamma:=\{\chi_A:A\in \Sigma\} \sub B_{L_\infty(\mu)}=B_{L_1(\mu)^*}$. 
Then ${\rm span}(\Gamma)$, which is the set of all simple functions,
is norm dense in $L_\infty(\mu)$. Condition~(ii) says that $(f_n)_{n\in \N}$ converges to~$0$ with respect to $\sigma(L_1(\mu),\Gamma)$.
Since $(f_n)_{n\in \N}$ is bounded (by condition~(i)), we can apply Lemma~\ref{lem:weakconvergence}
to conclude that $(f_n)_{n\in \N}$ is
weakly null in~$L_1(\mu)$. 

Finally, by the weak Banach--Saks property of~$L_1(\mu)$, there is a strictly increasing sequence $(n_j)_{j\in \N}$ in $\N$ such that
\[
 \lim_{m\to \infty} \left\|\frac1m\sum_{j=1}^{m}\xi_{n_j}\right\|_{\ca(\Sigma)}
 = 
 \lim_{m\to \infty} \left\|\frac1m\sum_{j=1}^{m}f_{n_j}\right\|_{L_1(\mu)}=0,
\]
as required.
\end{proof}

\begin{cor}\label{cor:weaknull}
Let $\nu:\Sigma\to c_0$ be a vector measure defined on a $\sigma$-algebra~$\Sigma$. 
For each $n\in \N$ we consider the composition $\nu_n:=e_n^*\circ \nu \in {\rm ca}(\Sigma)$. 
Then there is a strictly increasing sequence $(n_j)_{j\in \N}$ in $\N$ such that
\[
 \lim_{m\to \infty} \left\|\frac1m\sum_{j=1}^{m}\nu_{n_j}\right\|_{\ca(\Sigma)}
 = 0.
\]
\end{cor}
\begin{proof}
Since $\nu$ takes values in~$c_0$, we have
$$
	\lim_{n\to \infty}\nu_n(A)=0
	\quad
	\text{for every $A\in \Sigma$}.
$$
This condition is equivalent to the fact that $(\nu_n)_{n\in \N}$ is weakly null in ${\rm ca}(\Sigma)$
(Theorem~\ref{theo:ca}). The conclusion now follows from the weak Banach--Saks property
of ${\rm ca}(\Sigma)$ (Theorem~\ref{theo:Szlenk}).
\end{proof}

The idea of the proof of Theorem~\ref{theo:c0} is that
the weak Banach--Saks property ensures that the coordinate measures of a $c_0$-valued vector measure
have arithmetic means of arbitrarily small norm, 
while Schreier blocks force the norm of the same averages to be at least~$1$.

\begin{proof}[Proof of Theorem~\ref{theo:c0}]
Let $\Sigma$ be a $\sigma$-algebra and let $\nu:\Sigma \to c_0$ be a vector measure. We will prove that
$K_{\mathcal S}$ is not contained in $\mathcal{R}(\nu)$.

Let $(\nu_n)_{n\in \N}$ and $(n_j)_{j\in \N}$ be as in Corollary~\ref{cor:weaknull}. For each $m\in \N$ we write
$$
	s_m:=\sum_{j=1}^m \nu_{n_j} \in {\rm ca}(\Sigma),
$$ 
so that
$\lim_{m\to \infty} \|\frac1ms_m\|_{{\rm ca}(\Sigma)}=0$. Since for each $m\in \N$ we have
\[
 \left\|\frac1m\sum_{j=m+1}^{2m}\nu_{n_j}\right\|_{{\rm ca}(\Sigma)}
 =\left\|\frac{1}{m} (s_{2m}-s_m)\right\|_{{\rm ca}(\Sigma)}
 \le 2\left\|\frac{1}{2m}s_{2m}\right\|_{{\rm ca}(\Sigma)}+\left\|\frac{1}{m}s_m\right\|_{{\rm ca}(\Sigma)},
\]
we get 
\begin{equation}\label{eqn:limit}
	\lim_{m\to \infty} \left\|\frac1m\sum_{j=m+1}^{2m}\nu_{n_j}\right\|_{{\rm ca}(\Sigma)}=0.
\end{equation}

By contradiction, suppose that $K_{\mathcal S} \sub \mathcal{R}(\nu)$.
Fix $m\in\N$. The set
\[
 F_m:=\{n_{m+1},\ldots,n_{2m}\}
\]
belongs to~$\mathcal S$, because $|F_m|=m <n_{m+1}=\min F_m$. Therefore there is $A_m\in\Sigma$ such that
\[
  \nu(A_m)=\chi_{F_m}
\]
and so
\[
 \nu_{n_j}(A_m)=e_{n_j^*}(\nu(A_m))=e_{n_j^*}(\chi_{F_m})=1
 \quad
 \text{for every $j\in \{m+1,\dots,2m\}$.}
\]
It follows that
\begin{equation}\label{eqn:no-limit}
 \left\|\frac1m\sum_{j=m+1}^{2m}\nu_{n_j}\right\|_{{\rm ca}(\Sigma)}
 \ge
 \left|\left(\frac1m\sum_{j=m+1}^{2m}\nu_{n_j}\right)(A_m)\right|
 =\frac1m\sum_{j=m+1}^{2m}\nu_{n_j}(A_m)=1.
\end{equation}
Now, the validity of~\eqref{eqn:no-limit} for all $m\in \N$ contradicts \eqref{eqn:limit}.
\end{proof}

De Kock and Puglisi (see \cite[Theorem~3.2]{koc-pug}) showed that if a countable subset of $c_0$ lies in the range of a vector measure taking values in
a Banach space containing~$c_0$, then such a set is contained in the range of a $c_0$-valued vector measure. Their proof can be modified slightly to yield 
an improvement, see Theorem~\ref{theo:KP} below. 

Recall that a Banach space~$X$ is called \emph{weakly compactly generated} if there is a weakly compact set $K \sub X$
such that  $X=\overline{{\rm span}}(K)$. We will use the result of Amir and Lindenstrauss that every weakly compactly generated Banach space~$X$
has the \emph{separable complementation property}, that is, for every separable subspace $Y \sub X$ there is a projection on~$X$ 
whose range is separable and contains~$Y$ (see, e.g., \cite[Theorem~13.11 and Proposition~13.13]{fab-ultimo}).

\begin{theo}\label{theo:KP}
Let $X$ be a Banach space containing a subspace~$Y$ isomorphic to~$c_0$. Let $\nu :\Sigma\to X$ be a 
vector measure defined on a $\sigma$-algebra~$\Sigma$. Then there is a vector measure
$\tilde{\nu}:\Sigma\to Y$ such that
$$
	\tilde{\nu}(A)=\nu(A)  \qquad\text{whenever }A\in\Sigma\text{ and }\nu(A)\in Y.
$$
Consequently, $\mathcal{R}(\nu) \cap Y\sub \mathcal{R}(\tilde{\nu})$.
\end{theo}
\begin{proof}
Since $Y$ is separable, it is weakly compactly generated; let $K \sub Y$ be a weakly compact set such that $Y=\overline {{\rm span}}(K)$.
Let $L \sub X$ be the weak closure of~$\mathcal{R}(\nu)$. Then $L$ is weakly compact. Hence,
$L\cup K$ is weakly compact in~$X$ and therefore the subspace
\[
 	Z:=\overline{\operatorname{span}}\bigl( \mathcal{R}(\nu)\cup Y\bigr)=\overline{\operatorname{span}}\bigl(L \cup K\bigr)\sub X
\]
is weakly compactly generated. 

The Amir--Lindenstrauss theorem above yields a projection
$P$ on~$Z$ such that $Z_0:=P(Z)$ is separable and contains~$Y$. Note that $Y$ is a subspace of the separable Banach space~$Z_0$ and 
is isomorphic to~$c_0$.
By Sobczyk's theorem, there is a projection $Q$ on~$Z_0$
such that $Q(Z_0)=Y$ (see, e.g., \cite[Corollary~2.5.9]{alb-kal}). 

Let us consider the operator $T:=Q\circ P: Z \to Y$. Since $\nu$ takes values in~$Z$, we can define
a vector measure $\tilde{\nu}:\Sigma\to Y$ by 
$$
	\tilde{\nu}:=T\circ \nu=Q\circ P \circ \nu.
$$
Note that $\tilde{\nu}$ satisfies the required property: if $\nu(A)\in Y$ for some $A\in \Sigma$, then
\[
  \tilde{\nu}(A)=Q(P(\nu(A)))=Q(\nu(A))=\nu(A),
\]
because both $P$ and $Q$ coincide with the identity when restricted to~$Y$.
\end{proof}

\begin{theo}\label{theo:subspace}
If $X$ is a Banach space having property~AD, then $X$ contains no subspace isomorphic to~$c_0$. 
\end{theo}
\begin{proof}
By contradiction, suppose that $X$ has property~AD and contains a subspace~$Y$ isomorphic to~$c_0$. 
Let $T:c_0\to Y$ be an isomorphism and consider the weakly compact set $K:=T(K_{\mathcal S}) \sub Y$. 
Then there is an $X$-valued vector measure~$\nu$ such that $K \sub \mathcal{R}(\nu)$. Now,
Theorem~\ref{theo:KP} ensures the existence of a $Y$-valued vector measure~$\tilde{\nu}$ such that 
$K \sub \mathcal{R}(\tilde{\nu})$. Then $T^{-1}\circ \tilde{\nu}$ is a $c_0$-valued vector measure satisfying $K_{\mathcal{S}}\sub \mathcal{R}(T^{-1}\circ \tilde{\nu})$,
which contradicts Theorem~\ref{theo:c0}.  
\end{proof}

\begin{pro}\label{pro:dualc0}
If $X$ is a Banach space having property~AD, then $X^*$ contains no subspace isomorphic to~$c_0$.
\end{pro}
\begin{proof}
On the one hand, for an arbitrary Banach space~$X$, its dual $X^*$ contains no isomorphic copy of~$c_0$
if and only if $X$ contains no complemented isomorphic copy of~$\ell_1$ (see, e.g., \cite[Theorem~4.68]{ali-bur}).
On the other hand, as we mentioned in the introduction: (i)~property~AD is inherited by complemented subspaces, and (ii)~the space $\ell_1$ fails it. 
Hence, if $X$ has property~AD, then it contains no complemented isomorphic copy of~$\ell_1$ and the conclusion follows. 
\end{proof}

\begin{cor}\label{cor:lattice}
Let $X$ be a Banach lattice. If $X$ has property~AD, then it is superreflexive. 
\end{cor}
\begin{proof}
It is known that every reflexive Banach space having property~AD is superreflexive (see the introduction). It therefore remains to prove
that a Banach lattice having property~AD is reflexive. This is a consequence of Theorem~\ref{theo:subspace} and Proposition~\ref{pro:dualc0}, because
a Banach lattice $X$ is reflexive if and only if both $X$ and $X^*$ contain no isomorphic copy of~$c_0$
(see, e.g., \cite[Theorems 4.61 and~4.70]{ali-bur}).
\end{proof}

\subsection*{Statement on AI use}
During the preparation of this paper, the author used Microsoft Copilot, powered by GPT-5.6 Sol, to explore possible proof strategies and to improve the exposition. The author independently reconstructed and verified every proof, checked every cited result against its original source, and determined the final mathematical content and presentation. The author takes full responsibility for the originality, correctness, and integrity of the paper.

\bibliographystyle{amsplain}

\end{document}